
\newcount\secno
\newcount\prmno
\def\section#1{\vskip1truecm
               \global\def\currenvir{section}
               \global\advance\secno by1\global\prmno=0
               {\bf \number\secno. {#1}}
               \smallskip}

\def\subsection{\global\def\currenvir{subsection}
                \global\advance\prmno by1
               \smallskip  \ind{ (\number\secno.\number\prmno) }}
\def\subsec{\global\def\currenvir{subsection}
                \global\advance\prmno by1\smallskip
                { (\number\secno.\number\prmno)\ }}

\def\proclaim#1{\global\advance\prmno by 1
                {\bf #1 \the\secno.\the\prmno$.-$ }}

\long\def\th#1 \enonce#2\endth{%
   \medbreak\proclaim{#1}{\it #2}\global\def\currenvir{th}\smallskip}

\def\rem#1{\global\advance\prmno by 1
{\it #1} \the\secno.\the\prmno$.-$ }

\magnification 1250 \pretolerance=500 \tolerance=1000
\brokenpenalty=5000 \mathcode`A="7041 \mathcode`B="7042
\mathcode`C="7043 \mathcode`D="7044 \mathcode`E="7045
\mathcode`F="7046 \mathcode`G="7047 \mathcode`H="7048
\mathcode`I="7049 \mathcode`J="704A \mathcode`K="704B
\mathcode`L="704C \mathcode`M="704D \mathcode`N="704E
\mathcode`O="704F \mathcode`P="7050 \mathcode`Q="7051
\mathcode`R="7052 \mathcode`S="7053 \mathcode`T="7054
\mathcode`U="7055 \mathcode`V="7056 \mathcode`W="7057
\mathcode`X="7058 \mathcode`Y="7059 \mathcode`Z="705A
\def\spacedmath#1{\def\packedmath##1${\bgroup\mathsurround =0pt##1\egroup$}
\mathsurround#1
\everymath={\packedmath}\everydisplay={\mathsurround=0pt}}
 \spacedmath{2pt}

\def\iso{\vbox{\hbox to .8cm{\hfill{$\scriptstyle\sim$}\hfill}
\nointerlineskip\hbox to .8cm{{\hfill$\longrightarrow $\hfill}} }}
\def\sdir_#1^#2{\mathrel{\mathop{\kern0pt\oplus}\limits_{#1}^{#2}}}

\font\eightrm=cmr8 \font\sixrm=cmr6

\def\pc#1{\tenrm#1\sevenrm}
\def\tx{\kern-1.5pt -}
\def\cqfd{\kern 2truemm\unskip\penalty 500\vrule height 4pt depth 0pt width
4pt\medbreak} 
\def\no{n\up{o}\kern 2pt}
\def\ind{\par\hskip 1truecm\relax}

\font\pal=cmsy7

\def\sp#1{{\cal S}\kern-1pt\raise-1pt\hbox{\pal P}^{}_C(#1)}

\frenchspacing
\input xy
\xyoption{all}
\input amssym.def
\input amssym
\vsize = 25truecm \hsize = 16.1truecm \voffset = -.5truecm
\parindent=0cm
\baselineskip15pt \overfullrule=0pt

\vglue 2.5truecm \font\Bbb=msbm10

\centerline{\bf On the scope of validity of the norm}
\smallskip \centerline{\bf limitation theorem for quasilocal
fields}
\bigskip

\centerline{I.D. Chipchakov\footnote{$^{\ast}$}{Partially
supported by Grant MI-1503/2005 of the Bulgarian Foundation for
Scientific Research.}}
\par
\vskip1.truecm
\centerline{{\bf 1. Introduction}}
\par
\medskip
This paper is devoted to the study of norm groups of the fields
pointed out in the title, i.e. of fields whose finite extensions
are primarily quasilocal (briefly, PQL). It concentrates on the
special case where the considered ground fields are strictly
quasilocal, i.e. their finite extensions are strictly PQL (or
equivalently, these extensions admit one-dimensional local class
field theory, see [7]). The paper shows (see Theorem 1.1) that if
$E$ is a quasilocal field, $R/E$ is a finite separable extension,
and $R _{\rm ab}$ is the maximal abelian subextension of $E$ in
$R$, then the norm groups $N(R/E)$ and $N(R _{\rm ab}/E)$ are
equal, provided that the natural Brauer group homomorphism Br$(E)
\to $ Br$(L)$ is surjective, for every finite extension $L$ of
$E$. This is established in a more general form used in [10] (see
also (5.2) (i)) for describing the norm groups of finite
separable extensions of strictly quasilocal fields with Henselian
discrete valuations. Relying on [10], we prove here that Theorem
1.1 and the main results of [9], stated as (1.1) (ii), determine
to a considerable extent the scope of validity of the classical
norm limitation theorem (cf. [11, Ch. 6, Theorem 8]), in the case
of strictly PQL ground fields. The present research also sheds
light on the possibility of reducing the study of norm groups of
quasilocal fields to the special case of finite abelian
extensions.
\par
\medskip
The basic field-theoretic notions needed for describing the main
results of this paper are the same as those in [9]. As usual, $E
^{\ast }$ denotes the multiplicative group of a field $E$. We say
that $E$ is formally real, if $-1$ is not presentable as a finite
sum of squares of elements of $E$; the field $E$ is called
nonreal, otherwise. For convenience of the reader, we recall that
$E$ is said to be a PQL-field, if every cyclic extension $F$ of
$E$ is embeddable as an $E$-subalgebra in each central division
$E$-algebra $D$ of Schur index ind$(D)$ divisible by the degree
$[F\colon E]$. When this occurs, we say that $E$ is strictly PQL, if the
$p$-component Br$(E) _{p}$ of the Brauer group Br$(E)$ is
nontrivial in case $p$ runs through the set $P(E)$ of those prime
numbers, for which $E$ is properly included in its maximal
$p$-extension $E (p)$ in a separable closure $E _{\rm sep}$ of
$E$. It is worth noting that PQL-fields and quasilocal fields
appear naturally in the process of characterizing some of the
basic types of stable fields with Henselian valuations (see [8]
and the references there). Our research, however, is primarily
motivated by the fact that strictly PQL-fields admit local class
field theory, and by the validity of the converse in all presently
known cases (see [7, Theorem 1 and Sect. 2]). As to the choice of
our main topic, it is determined to a considerable extent by the
following results:
\par
\medskip
(1.1) (i) $N(R/E) = N(R _{\rm ab}/E)$, provided that $R$ is a
finite separable extension of a field $E$ possessing a Henselian
discrete valuation with a quasifinite residue field $\widehat E$
[18] (see also [24] and [29]]);
\par
(ii) $N(R/E) = N(R _{\rm ab}/E)$ in case $E$ is a PQL-field and
$R$ is an intermediate field of a finite Galois extension $M/E$
with a nilpotent Galois group; for each nonnilpotent finite group
$G$, there exists an algebraic extension $E(G)$ of $\hbox{\Bbb
Q}$, which is strictly PQL and has a Galois extension $M(G)$, such
that $G(M(G)/E(G))$ is isomorphic to $G$ and $N(M(G)/E(G))$ is a
proper subgroup of $N(M(G) _{\rm ab}/E(G))$ [9, Theorems 1.1 and
1.2];
\par
(iii) If $E$ is an algebraic strictly PQL-extension of a global
field $E _{0}$, and $R/E$ is a finite extension, then $N(R/E) =
N(\Phi (R)/E)$, for some abelian finite extension $\Phi (R)$ of
$E$, which is uniquely determined by $R/E$, up-to an
$E$-isomorphism (see the references after the statement of [9,
Theorem 1.2]).
\par
\medskip
The main purpose of this paper is to shed an additional light on
these facts by proving the following two statements (the former 
of which generalizes (1.1) (i), see also Remark 4.4, for more 
details):
\par
\medskip
{\bf Theorem 1.1.} {\it Let $E$ be a quasilocal field and $R/E$ a
finite separable extension. Then $N(R/E) = N(R _{\rm ab}/E)$ in
the following two special cases:}
\par
(i) {\it The natural homomorphism of {\rm Br}$(E)$ into {\rm
Br}$(L)$ is surjective, for every finite extension $L$ of $E$;}
\par
(ii) {\it There exists an abelian finite extension $\Phi (R)$ of
$E$, such that \par \noindent $N(\Phi (R)/E) = N(R/E)$.}
\par
\medskip
{\bf Theorem 1.2.} {\it There exists a strictly quasilocal nonreal
field $E$ satisfying the following conditions:}
\par
(i) {\it the absolute Galois group $G _{K} := G(K _{\rm sep}/K)$
is not pronilpotent;}
\par
(ii) {\it every finite extension $R$ of $K$ is subject to the
following alternative:}
\par
($\alpha $) $R$ {\it is an intermediate field of a finite Galois
extension $M(R)/K$ with a nilpotent Galois group;}
\par
($\beta $) $N(R/K)$ {\it does not equal the norm group of any
abelian finite extension of $K$.}
\par
\medskip
In addition to (1.1) and Theorems 1.1 and 1.2, it has been proved
in [6] that the description of the norm groups of finite 
separable extensions of a strictly PQL-field $F$ does not reduce 
to the study of Galois extensions $M/F$ with $G(M/F)$ belonging 
to any given proper class of finite groups, which is closed under 
the formation of subgroups, quotient groups and group extensions. 
Also, it has been shown in [5] that a formally real strictly 
quasilocal field $E$ has the properties required by Theorem 1.2 
(i) and (ii) unless it is real closed.
\par
\medskip
Throughout the paper, simple algebras are supposed to be
associative with a unit and finite-dimensional over their centres,
and Galois groups are viewed as profinite with respect to the
Krull topology. For each simple algebra $A$, we consider only
subalgebras of $A$ containing its unit. Our basic terminology and
notation concerning valuation theory, simple algebras and Brauer
groups are standard (for example, as in [12; 15; 36] and [20], as
well as those related to profinite groups, Galois cohomology,
field extensions and Galois theory (see, for example, [25; 13] and
[15]). We refer the reader to [28, Sect. 1] and [4, Sect. 2], for
the definitions of a symbol algebra and of a symbol $p$-algebra
(see also [26, Ch. XIV, Sects. 2 and 5]).
\par
\medskip
Here is an overview of the paper: Section 2 includes preliminaries
used in the sequel. Theorems 1.1 and 1.2 are proved in Sections
3-4 and 5, respectively. Section 5 contains a characterization of 
the fields singled out by Theorem 1.2 among those endowed with a 
Henselian discrete valuation and possessing the strictly 
PQL-property.
\par
\vskip 0.6 truecm \centerline {\bf 2. Preliminaries}
\par
\medskip
Let $E$ be a field, Nr$(E)$ the set of norm groups of finite
extensions of $E$ in $E _{\rm sep}$, and $\Omega (E)$ the set of
finite abelian extensions of $E$ in $E _{\rm sep}$. We say that
$E$ admits (one-dimensional) local class field theory, if the
mapping $\pi $ of $\Omega (E)$ into Nr$(E)$ defined by the rule
$\pi (F) = N(F/E)\colon \ F \in \Omega (E)$, is injective and 
satisfies the following two conditions, for each pair $(M _{1}, 
M _{2}) \in \Omega (E) \times \Omega (E)\colon $
\par
The norm group of the compositum $M _{1}M _{2}$ is equal to the
intersection $N(M _{1}/E) \ \cap $
\par \noindent
$N(M _{2}/E)$ and $N((M _{1} \cap M _{2})/E)$ equals the inner 
group product $N(M _{1}/E)N(M _{2}/E)$.
\par
We call $E$ a field with (one-dimensional) local $p$-class field
theory, for some prime $p$, if the restriction of $\pi $ on the
set of finite abelian extensions of $E$ in $E (p)$ has the same
properties. Our approach to the study of fields with such a theory
is based on the following two lemmas (proved in [9]).
\par
\medskip
{\bf Lemma 2.1.} {\it Let $E$ be a field and $L$ an extension of
$E$ presentable as a compositum of extensions $L _{1}$ and $L
_{2}$ of $E$ of relatively prime degrees. Then $N(L/E) = N(L
_{1}/E) \cap N(L _{2}/E)$, $N(L _{1}/E) = E ^{\ast } \cap N(L/L
_{2})$, and there is a group isomorphism $E ^{\ast }/N(L/E) \cong
(E ^{\ast }/N(L _{1}/E)) \times (E ^{\ast }/N(L _{2}/E))$.}
\par
\medskip
{\bf Lemma 2.2.} {\it Let $E$ be a field, $M$ a finite Galois
extension of $E$ with a nilpotent Galois group $G(M/E)$, $R$ an
intermediate field of $M/E$ not equal to $E$, $\Pi $ the set of
prime numbers dividing $[R\colon E]$, $M _{p}$ the maximal 
$p$-extension of $E$ in $M$, and $R _{p} = R \cap M _{p}$, for 
each $p \in \Pi $. Then:}
\par
(i) $R$ {\it is equal to the compositum of the fields $R 
_{p}\colon \ p \in \Pi $, and $[R\colon E] = \prod _{p \in \Pi } 
[R _{p}\colon E]$;}
\par
(ii) $N(R/E) = \cap _{p \in \Pi } N(R _{p}/E)$ {\it and the 
quotient group $E ^{\ast }/N(R/E)$ is isomorphic to the direct 
product of the groups $E ^{\ast }/N(R _{p}/E)\colon \ p \in \Pi $.}
\par
\medskip
It is clear from Lemma 2.2 that a field $E$ admits local class
field theory if and only if it admits local $p$-class field
theory, for every $p \in P(E)$. Our next lemma, proved in [8,
Sect. 4], shows that Br$(E) _{p} \neq \{0\}$ whenever $E$ is a
field with such a theory, for a given $p \in P(E)$.
\par
\medskip
{\bf Lemma 2.3.} {\it Let $E$ be a field, such that} Br$(E) _{p} =
\{0\}${\it , for some prime number p. Then} Br$(E _{1}) _{p} =
\{0\}$ {\it and $N(E _{1}/E) = E ^{\ast }$, for every finite
extension $E _{1}$ of $E$ in $E (p)$.}
\par
\medskip
The following lemma is known (cf. [25, Ch. II, 2.3 and 3.1]) and
plays an essential role in the proof of Theorem 1.1. 
\par
\medskip
{\bf Lemma 2.4.} {\it For a field $E$ and a prime number $p$, the
following conditions are equivalent:}
\par
(i) Br$(E ^{\prime }) _{p} = \{0\}${\it , for every algebraic
extension $E ^{\prime }$ of $E$;}
\par
(ii) {\it The exponent of the group $E _{1} ^{\ast }/N(E _{2}/E
_{1})$ is not divisible by $p$, for any pair $(E _{1}, E _{2})$ of
finite extensions of $E$ in $E _{\rm sep}$, such that $E _{1}
\subseteq E _{2}$.}
\par
\medskip
For a detailed proof of Lemma 2.4, we refer the reader to [5].
Let now $\Phi $ be a field and $\Phi _{p}$ the extension of $\Phi
$ in $\Phi _{\rm sep}$ generated by a primitive $p$-th root of
unity $\varepsilon _{p}$, for some prime $p$. It is well-known 
(cf. [15, Ch. VIII, Sect. 3]) that then $\Phi _{p}/\Phi $ is a 
cyclic extension of degree $[\Phi _{p}\colon \Phi ] := m$ dividing 
$p - 1$. Denote by $\varphi $ some $\Phi $-automorphism of $\Phi 
_{p}$ of order $m$, fix an integer $s$ so that $\varphi 
(\varepsilon _{p}) = \varepsilon _{p} ^{s}$, and put $V  _{i} =
\{\alpha _{i} \in \Phi _{p} ^{\ast }\colon \ \varphi (\alpha _{i})\alpha
_{i} ^{-s ^{i}} \in \Phi _{p} ^{\ast p}\}$ and $\overline V _{i}
= V _{i}/\Phi _{p} ^{\ast p}\colon \ i = 0, \dots , m - 1$. Clearly, the
quotient group $\Phi _{p} ^{\ast }/\Phi _{p} ^{\ast p} :=
\overline \Phi _{p} $ can be viewed as a vector space over the
field $\hbox{\Bbb F} _{p}$ with $p$ elements. Considering the
linear operator $\bar \varphi $ of $\overline \Phi _{p}$, defined
by the rule $\bar \varphi (\alpha \Phi _{p} ^{\ast p}) = \varphi
(\alpha )\Phi _{p} ^{\ast p}\colon \ \alpha \in \Phi _{p} ^{\ast }$,
and taking into account that the subspace of $\overline \Phi
_{p}$, spanned by its elements $\bar \varphi ^{i} (\bar \alpha 
)\colon \ i = 0, \dots , m - 1$, is finite-dimensional and $\bar 
\varphi $-invariant, for each $\bar \alpha \in \overline \Phi 
_{p}$, one obtains from Maschke's theorem the following statement:
\par
\medskip
(2.1) The sum of the subspaces $\overline V _{i}\colon \ i = 0, 
\dots , m - 1$ is direct and equal to $\overline \Phi _{p}$.
\par
\medskip
Let $L$ be an extension of $\Phi _{p}$ in $\Phi _{\rm sep}$,
obtained by adjoining a $p$-th root $\eta _{p}$ of an element
$\beta \in (\Phi _{p} ^{\ast } \setminus \Phi _{p} ^{\ast p})$. It
is clear from Kummer's theory that $[L\colon \Phi ] = pm$ and
the following assertions hold true:
\par
\medskip
(2.2) $L/\Phi $ is a Galois extension if and only if $\beta \in V
_{j}$, for some index $j$. Such being the case, every $\Phi
_{p}$-automorphism $\psi $ of $L$ of order $p$ satisfies the
equality $\varphi ^{\prime }\psi \varphi ^{\prime -1} = \psi
^{s'}$, where $s ^{\prime } = s ^{1-j}$ and $\varphi ^{\prime }$
is an arbitrary automorphism of $L$ extending $\varphi $.
Moreover, $L$ and $\Phi $ are related as follows: 
\par
(i) $L/\Phi $ is cyclic if and only if $\beta \in V _{1}$ (Albert
[1, Ch. IX, Theorem 6]);
\par
(ii) $L$ is a root field over $\Phi $ of the binomial $X ^{p} -
a$, for some $a \in \Phi ^{\ast }$, if and only if $\beta \in V
_{0}$, i.e. $s ^{\prime } = s$; when this occurs, one can take as
$a$ the norm $N _{\Phi } ^{\Phi _{p}} (\beta )$.
\par
\medskip
Statements (2.1), (2.2) and the following observations will be
used for proving Theorems 1.1 and 1.2.
\par
\medskip
(2.3) For a symbol $\Phi _{p}$-algebra $A _{\varepsilon _{p}}
(\alpha , \beta ; \Phi _{p})$ (of dimension $p ^{2}$), where
$\alpha \in \Phi _{p} ^{\ast }$ and $\beta \in V _{j} \setminus
\Phi _{p} ^{\ast p}$, the following conditions are equivalent:
\par
(i) $A _{\varepsilon _{p}} (\alpha , \beta ; \Phi _{p})$ is $\Phi
_{p}$-isomorphic to $D \otimes _{\Phi } \Phi _{p}$, for some
central simple $\Phi $-algebra $D$;
\par
(ii) If $\alpha = \prod _{i=0} ^{m-1} \alpha _{i}$ and $\alpha
_{i} \in V _{i}$, for each index $i$, then $A _{\varepsilon _{p}}
(\alpha , \beta ; \Phi _{p})$ is isomorphic to the symbol $\Phi
_{p}$-algebra $A _{\varepsilon _{p}} (\alpha _{j'}, \beta ; \Phi
_{p})$, where $j ^{\prime }$ is determined so that $m$ divides $j
^{\prime } + j - 1$;
\par
(iii) With notations being as in (ii), $\alpha _{i} \in N(L/\Phi
_{p})\colon \ i \neq j ^{\prime }$.
\par
\medskip
The main results of [7, Sect. 2] and [8] used in the present paper
(sometimes without an explicit reference) can be stated as follows:
\par
\medskip
{\bf Proposition 2.5.} {\it Let $E$ be a strictly $p$-quasilocal
field, for some $p \in P(E)$. Assume also that $R$ is a finite
extension of $E$ in $E (p)$, and $D$ is a central division
$E$-algebra of $p$-primary dimension. Then $R$, $E$ and $D$ have 
the following properties:}
\par
(i) {\it $R$ is a $p$-quasilocal field and  ${\rm ind}(D) = {\rm 
exp}(D)$;}
\par
(ii) Br$(R) _{p}$ {\it is a divisible group unless $p = 2$, $R = E$
and $E$ is formally real; in the noted exceptional case, {\rm
Br}$(E) _{2}$ is of order $2$ and $E (2) = E(\sqrt{-1})$;}
\par
(iii) $E$ {\it admits local $p$-class field theory, provided that
Br$(E) _{p} \neq \{0\}$;}
\par
(iv) {\it $R$ embeds in $D$ as an $E$-subalgebra if and only if
$[R\colon E]$ divides {\rm ind}$(D)$.}
\par
\vskip2.truecm \centerline {\bf 3. $p$-primary analogue to
Theorem 1.1}
\par
\medskip
Let $E$ be a field, $R/E$ a finite separable extension, and for
each prime $p$, let $R _{{\rm ab},p}$ be the maximal abelian
$p$-extension of $E$ in $R$, $\rho _{p}$ the greatest integer
dividing $[R\colon E]$ and not divisible by $p$, and $N _{p} (R/E)$ the
set of those elements $u _{p} \in E ^{\ast }$, for which the
co-set $u _{p}N(R/E)$ is a $p$-element of the  group $E ^{\ast
}/N(R/E)$. Clearly, $u ^{\rho _{p}} \in N _{p} (R/E)$, for every
$u \in E ^{\ast }$. Observing also that $u ^{\rho _{p}} \in N(R
_{{\rm ab},p}/E)$ whenever $u \in N(R _{\rm ab}/E)$ and $p$ is
prime, one concludes that Theorem 1.1 (i) can be deduced from its
$p$-primary analogue stated as follows:
\par
\medskip
{\bf Theorem 3.1.} {\it Assume that $E$ is a quasilocal field,
such that the natural homomorphism of} Br$(E)$ {\it into} Br$(L)$
{\it maps} Br$(E) _{p}$ {\it surjectively on} Br$(L) _{p}${\it ,
for some prime number $p$ and every finite extension $L$ of $E$.
Then $N(R/E) = N(R _{{\rm ab},p}/E) \cap $
\par \noindent 
$N _{p} (R/E)$, for each finite extension $R$ of $E$ in $E 
_{\rm sep}$.}
\par
\medskip
In what follows, up-to the end of the next Section, our main
objective is to prove Theorems 3.1 and 1.1. Evidently, $N(R/E)
\subseteq (N _{p} (R/E)$ $\cap N(R _{\rm ab}/E))$, so we have to
prove that $N _{p} (R/E) \cap N(R _{\rm ab}/E)$ is a subgroup of
$N(R/E)$. Our assumptions show that if Br$(E) _{p} = \{0\}$, then
Br$(L) _{p} = \{0\}$, for every finite extension $L$ of $E$, which
reduces our assertion to a consequence of Lemma 2.4. Assuming
further that Br$(E) _{p} \neq \{0\}$ and $\hbox{\Bbb F} _{p}$ is
a field with $p$ elements (identifying it with the prime subfield
of $E$, in the case of char$(E) = p$), we prove in the rest of
this Section the validity of Theorem 3.1 in the special case where
$R/E$ is a normal extensions with a solvable Galois group. The
main part of our argument is presented by the following two
lemmas.
\par
\medskip
{\bf Lemma 3.2.} {\it Let $E$ be a field and $p$ a prime number
satisfying the conditions of Theorem 3.1, and let $M/E$ be a
Galois extension with $G(M/E)$ satisfying the following 
conditions:}
\par
(i) $G(M/E)$ {\it is nonabelian and isomorphic to a semidirect
product $E _{p;k} \times C _{\pi }$ of an elementary abelian
$p$-group of order $p ^{k}$ by a group $C _{\pi }$ of prime order
$\pi $ not equal to $p$, where $k$ is the minimal positive 
integer solution to the congruence $p ^{k} \equiv 1$ $({\rm mod} 
\pi )$;}
\par
(ii) $E _{p;k}$ {\it is a minimal normal subgroup of $G(M/E)$.
\par
Then $N(M/E _{1})$ includes $E ^{\ast }$, where $E _{1}$ is the 
intermediate field of $M/E$ corresponding by Galois theory to $E 
_{p;k}$.}
\par
\medskip
{\it Proof.} Our assumptions indicate that $E _{1}/E$ is a cyclic
extension of degree $\pi $, and under the additional hypothesis
that Br$(E) _{p} \neq \{0\}$, this means that Br$(E _{1}) _{p}
\neq \{0\}$ (see [20, Sect. 13.4]). Therefore, by Proposition 2.5
(iii), $E _{1}$ admits local $p$-class field theory, so it is
sufficient to show that $E ^{\ast } \subseteq N(M _{1}/E _{1})$,
for every cyclic extension $M _{1}$ of $E _{1}$ in $M$. Suppose
first that $E$ contains a primitive $p$-th root of unity or
char$(E) = p$, and fix an $E$-automorphism $\psi $ of $E _{1}$ of
order $\pi $. As $G(M/E _{1})$ is an elementary abelian $p$-group
of rank $k$, Kummer$^{,}$s theory and the Artin-Schreier theorem
(cf. [15, Ch. VIII, Sect. 6]) imply the existence of a subset $S =
\{\rho _{j}\colon \ j = 1, \dots , k\}$ of $E _{1}$, such that the root field
over $E _{1}$ of the polynomial set $\{f _{j} (X) = X ^{p} - uX -
\rho _{j}\colon \ j = 1, \dots , k\}$ equals $M$, where $u = 1$, if char$(E)
= p$, and $u = 0$, otherwise. For each index $j$, denote by $z
_{j}$ the element $\psi (u _{j})u _{j} ^{-1}$ in case $E$ contains
a primitive $p$-th root of unity, and put $z _{j} = \psi (u _{j})
- u _{j}$, if char$(E) = p$. Note that $M$ is a root field over $E
_{1}$ of the set of polynomials $\{g _{j} (X) = X ^{p} - uX - z
_{j}\colon \ j = 1, \dots , k\}$. This can be deduced from the following two
statements:
\par
\medskip
(3.1) (i) If char$(E) = p$, $r(E _{1}) = \{\lambda ^{p} - \lambda
\colon \ \lambda \in E _{1}\}$ and $M(E _{1})$ is the additive subgroup
of $E$ generated by the union $S \cup r(E _{1})$, then $r(E _{1})$
and $M(E _{1})$ are $\psi $-invariant, regarded as vector spaces
over $\hbox{\Bbb F} _{p}$; moreover, the linear operator of the
quotient space $M(E _{1})/r(E _{1})$, induced by $\psi - id _{E
_{1}}$ is an isomorphism;
\par
(ii) If $E$ contains a primitive $p$-th root of unity, $M(E _{1})$
is the multiplicative subgroup of $E _{1} ^{\ast }$ generated by
the union $S \cup E _{1} ^{\ast p}$, and the mapping $\psi 
_{1}\colon \ E _{1} ^{\ast }/E _{1} ^{\ast p} \to E _{1} ^{\ast 
}/E _{1} ^{\ast p}$ is defined by the rule $\psi _{1} (\alpha  E 
_{1} ^{\ast p}) = \psi (\alpha )\alpha ^{-1}E _{1} ^{\ast p}\colon 
\ \alpha \in E _{1} ^{\ast }$, then $\psi _{1}$ is a linear 
operator of $E _{1} ^{\ast }/E _{1} ^{\ast p}$ (regarded as a 
vector space over $\hbox{\Bbb F} _{p}$), $M(E _{1})/E _{1} ^{\ast 
p}$ is a $k$-dimensional $\psi _{1}$-invariant subspace of $E _{1} 
^{\ast }/E _{1} ^{\ast p}$, and the linear operator of $M(E 
_{1})/E _{1} ^{\ast p}$ induced by $\psi _{1}$ is an isomorphism.
\par
\medskip
Most of the assertions of (3.1) are well-known. One should,
possibly, only note here that the concluding parts of (3.1) (i)
and (3.1) (ii) follow from the fact that $G(M/E _{1})$ is the
unique normal proper subgroup of $G(M/E)$, and by Galois theory,
this means that $E _{1}$ is the unique normal proper extension of
$E$ in $M$. The obtained result implies the nonexistence of a
cyclic extension of $E$ in $M$ of degree $p$, which enables one to
deduce from Kummer's theory and the Artin-Schreier theorem the
triviality of the kernels of the considered linear operators. Thus
our argument leads to the conclusion that the discussed special
case of Lemma 3.2 will be proved, if we establish the validity of
the following two statements, for each index $j$:
\par
\medskip
(3.2) (i) If $E$ contains a primitive $p$-th root of unity
$\varepsilon $ and $c$ is an element of $E ^{\ast }$, then the
symbol $E _{1}$-algebra $A _{\varepsilon } (z _{j}, c; E _{1})$ is
trivial;
\par
(ii) If char$(E) = p$ and $c \in E ^{\ast }$, then the $p$-symbol
$E$-algebra $E[z _{j}, c)$ is trivial.
\par
\medskip
Denote by $D _{j}$ the symbol $p$-algebra $E _{1} [\rho _{j}, c)$,
if char$(E) = p$, and the symbol $E _{1}$-algebra $A _{\varepsilon
} (\rho _{j}, c; E _{1})$ in case $E _{1}$ contains a primitive
$p$-th root of unity $\varepsilon $. It follows from the
assumptions of Theorem 3.1 that $D _{j}$ is isomorphic over $E
_{1}$ to $\Delta _{j} \otimes _{E} E _{1}$, for some central
division $E$-algebra $\Delta _{j}$. This implies that $\psi $ is
extendable to an automorphism $\bar \psi $ of $D _{j}$, regarded
as an algebra over $E$. Thus it becomes clear that $D _{j}$ is $E
_{1}$-isomorphic to $E _{1} [\psi (\rho _{j}), c)$ or $A
_{\varepsilon } (\psi (\rho _{j}), c; E _{1})$ depending on
whether or not char$(E) = p$. Applying now the general properties
of local symbols (cf. [26, Ch. XIV, Propositions 4 and 11]), one
proves (3.2).
\par
It remains for us to prove Lemma 3.2, assuming that $p \neq $
char$(E)$ and $E$ does not contain a primitive $p$-th root of
unity. Let $\varepsilon $ be such a root in $M _{\rm sep}$. It is
easily verified that if $E(\varepsilon ) \cap E _{1} = E$, then
$M(\varepsilon )/E(\varepsilon )$ is a Galois extension, such that
$G((M(\varepsilon )/E(\varepsilon ))$ is canonically isomorphic to
$G(M/E)$. Since $E(\varepsilon )$ and $p$ satisfy the conditions
of the lemma, our considerations prove in this case that
$E(\varepsilon ) ^{\ast } \subseteq N(M(\varepsilon )/E _{1}
(\varepsilon ))$. Hence, by Lemma 2.1, applied to the triple $(E
_{1}, M, E _{1} (\varepsilon ))$ instead of $(E, L _{1}, L _{2})$,
we have $E ^{\ast } \subseteq N(M/E _{1})$, which reduces the
proof of Lemma 3.2 to the special case in which $E _{1}$ is an
intermediate field of $E(\varepsilon )/E$. Fix a generator
$\varphi $ of $G(E(\varepsilon )/E)$, and an integer $s$ so that 
$\varphi (\varepsilon ) = \varepsilon ^{s}$. Observing that
$M/E$ is a noncyclic Galois extension of degree $p\pi $, one
obtains from (2.2) and the cyclicity of $M$ over $E _{1}$ that
$M(\varepsilon )$ is generated over $E(\varepsilon )$ by a $p$-th
root of an element $\rho $ of $E(\varepsilon )$ with the property
that $\varphi (\rho )\rho ^{-s'} \in E(\varepsilon ) ^{\ast p}$,
where $s ^{\prime }$ is a positive integer such that $s ^{\prime
\pi } \equiv s ^{\pi }$(mod $p$) and $s ^{\prime } \not\equiv
s$(mod $p$). It is therefore clear from (2.3), the surjectivity of
the natural homomorphism of Br$(E) _{p}$ into Br$(E(\varepsilon ))
_{p}$, and [20, Sect. 15.1, Proposition b] that $A _{\varepsilon }
(\rho , c; E(\varepsilon ))$ is isomorphic to the matrix
$E(\varepsilon )$-algebra $M _{p} (E(\varepsilon ))$, for every $c
\in E ^{\ast }$. One also sees that $E ^{\ast } \subseteq
N(M(\varepsilon )/E(\varepsilon ))$. As $[M\colon E _{1}] = p$ and
$[E(\varepsilon )\colon E _{1}]$ divides $(p - 1)/\pi $, Lemma 2.1
ensures now that $E ^{\ast } \subseteq N(M/E _{1})$, so Lemma 3.2
is proved.
\par
\medskip
{\bf Lemma 3.3.} {\it Assume that $E$ is a quasilocal field whose
finite extensions satisfy the conditions of Theorem 3.1, for
a given prime number $p$, and suppose that $M/E$ is a finite
Galois extension, such that $G(M/E)$ is a solvable group. Then $N
_{p} (M/E) \cap N(M _{{\rm ab},p}/E)$ is a subgroup of $N(M/E)$.}
\par
\medskip
{\it Proof.} It is clearly sufficient to prove the lemma under the
hypothesis that $N(M ^{\prime }/E ^{\prime })$ includes $N _{p} (M
^{\prime }/E ^{\prime }) \cap N(M ^{\prime } _{{\rm ab},p}/E
^{\prime })$, provided that $E ^{\prime }$ and $p$ satisfy the
conditions of Theorem 3.1, and $M ^{\prime }/E ^{\prime }$ is a
Galois extension with a solvable Galois group of order less than
$[M\colon E]$. As in the proof of [9, Theorem 1.1], we first show that
then one may assume further that $G(M/E)$ is a Miller-Moreno
group (i.e. nonabelian with abelian proper subgroups). Our 
argument relies on the fact that the class of fields satisfying 
the conditions of Theorem 3.1 is closed under the formation of 
finite extensions. Note that if $G(M/E)$ is not Miller-Moreno, 
then it possesses a nonabelian subgroup $H$ whose commutator 
subgroup $[H, H]$ is normal in $G(M/E)$. Indeed, one can take as 
$H$ the commutator subgroup $[G(M/E), G(M/E)]$ in case $G(M/E)$ 
is not metabelian, and suppose that $H$ is any nonabelian maximal 
subgroup of $G(M/E)$, otherwise. Denote by $F$ and $L$ the 
intermediate fields of $M/E$ corresponding to $H$ and $[H, H]$, 
respectively. Our choice of $H$ and Galois theory indicate that 
$L/E$ is a Galois extension such that $M _{\rm ab} \subseteq L$ 
and $E \neq L \neq M$, so our additional hypothesis and 
Lemma 2.2 lead to the conclusion that $N _{p} (L/E) \cap N(M 
_{{\rm ab},p}/E) = N _{p} (L/E) \cap N(M _{\rm ab}/E) \subseteq 
N(L/E)$ and $N _{p} (M/F) \cap N(L/F) \subseteq N(M/F)$. Let now 
$\mu $ be an element of $N _{p} (M/E) \cap N(M _{\rm ab,p}/E)$, 
and $\lambda \in L ^{\ast }$ a solution to the norm equation $N 
_{E} ^{L} (X) = \mu $. Then one can find an integer $k$ not 
divisible by $p$ and such that $N _{F} ^{L} (\lambda ) ^{k} \in 
N _{p} (M/F)$. It is therefore clear that $N _{F} ^{L} (\lambda ) 
^{k} \in N(M/F)$ and $\mu ^{k} \in N(M/E)$. As $\mu \in N _{p} 
(M/E)$, this implies that $\mu \in N(M/E)$, which yields the 
desired reduction. In view of the former part of (1.1) (ii),
one may also assume that $G(M/E)$ is a nonnilpotent Miller-Moreno
group. The assertion of Lemma 3.3 is obvious, if $p$ does not 
divide the order $o([G, G])$ of $[G, G]$, so we suppose further 
that $p \ \vert \ o([G, G])$. By the classification of these groups 
[17] (cf. also [22, Theorem 445]), this means that $G(M/E)$ has 
the following structure:
\par
\medskip
(3.3) (i) $G(M/E)$ is isomorphic to a semi-direct product $E
_{p;k} \times C _{\pi ^{n}}$ of $E _{p;k}$ by a cyclic group $C
_{\pi ^{n}}$ of order $\pi ^{n}$, for some different prime numbers
$p$ and $\pi $, where $k$ satisfies condition (i) of Lemma 3.2;
\par
(ii) $E _{p;k}$ is a minimal normal subgroup of $G(M/E)$, $E
_{p;k} = [G(M/E),$ $G(M/E)]$ and the centre of $G(M/E)$ equals the
subgroup $C _{\pi ^{n-1}}$ of $C _{\pi ^{n}}$ of order $\pi
^{n-1}$.
\par
\medskip
It follows from (3.2) and Galois theory that $M _{ab}/E$ is cyclic
of degree $\pi ^{n}$. This yields $N _{M _{ab}} ^{M} (\eta ) =
\eta ^{p ^{k}}$, for every $\eta \in M _{ab}$, and thereby,
implies that $c ^{p ^{k}} \in N(M/E)$ in case $c \in N(M _{\rm
ab}/E)$. It is therefore clear from the equality g.c.d.$(p ^{k},
\pi ^{n}) = 1$ that Lemma 3.3 will be proved, if we show that $c
^{\pi ^{n}} \in N(M/E)$ whenever $c \in E ^{\ast }$. By Lemma 3.2,
if $n = 1$, then $M ^{\ast }$ contains an element $\xi $ of norm
$c$ over $E _{1} = M _{ab}$, which means that $N _{E} ^{M} (\xi )
= c ^{\pi }$. Suppose now that $n \ge 2$, put $\tilde \pi = \pi
^{n-1}$, denote by $C _{\tilde \pi }$ the subgroup of $G(M/E)$ of
order $\tilde \pi $, and let $M ^{\prime }$ and $E ^{\prime }$ be
the intermediate fields of $M/E$ corresponding by Galois theory to
the subgroups $C _{\tilde \pi }$ and $E _{p;k}C _{\tilde \pi }$ of
$G(M/E)$, respectively. It is easily seen that $M ^{\prime }/E$ is
a Galois extension with $G(M ^{\prime }/E)$ satisfying the
conditions of Lemma 3.2, and $E ^{\prime }/E$ is a cyclic
extension of degree $\pi $. This ensures that $c ^{\pi } \in N(M
^{\prime }/E)$. Also, it becomes clear that $M = M ^{\prime }M
_{ab}$, $M ^{\prime } \cap M _{ab} = E ^{\prime }$ and $N _{M'}
^{M} (m ^{\prime }) = m ^{\prime \tilde \pi }$, for every $m
^{\prime } \in M ^{\prime }$. These observations show that $c
^{\pi ^{n}} \in N(M/E)$, so Lemma 3.3 is proved.
\par
\vskip0.6truecm \centerline{\bf 4. Proof of Theorems 3.1 and 1.1}
\par
\medskip
Retaining notation as in Section 3, we first consider the special
case in which $R$ is an intermediate field of a finite Galois
extension with a solvable Galois group. Our argument relies on
Lemma 3.3 and the following lemma.
\par
\medskip
{\bf Lemma 4.1.} {\it Under the hypotheses of Theorem 3.1, suppose
that $M/E$ is a Galois extension with a solvable Galois group
$G(M/E)$, and $R$ is an intermediate field of $M/E$, such that
$[R\colon E]$ is a power of $p$. Then $N(R/E) = N(R _{\rm ab}/E)$.}
\par
\medskip
{\it Proof.} Arguing by induction on $[M\colon E]$, one obtains from the
conditions of Theorem 3.1 that it is sufficient to prove the
lemma, assuming in addition that $N(R _{1}/E _{1}) = N(R ^{\prime
}/E _{1})$ whenever $E _{1}$ and $R _{1}$ are intermediate fields
of $M/E$, such that $E _{1} \neq E$, $E _{1} \subseteq R _{1}$,
$[R _{1}\colon E _{1}]$ is a power of $p$, and $R ^{\prime }$ is the
maximal abelian extension of $E _{1}$ in $R _{1}$. Suppose first
that $R _{\rm ab} \neq E$. Then the inductive hypothesis, applied
to the the pair $(E _{1}, R _{1}) = (R _{\rm ab}, R)$, gives
$N(R/E) = N(R ^{\prime }/E)$, and since $R ^{\prime }$ is a
subfield of the maximal $p$-extension $M _{p}$ of $E$ in $M$, this
enables one to obtain from the former part of (1.1) (ii) that $N(R
^{\prime }/E) = N(R _{\rm ab}/E)$.
\par
It remains to be seen that $N(R/E) = E ^{\ast }$ in the special
case of $R _{\rm ab} = E$. Our argument relies on the fact that $E
^{\ast }/N(M _{\rm ab}/E)$ is a group of exponent dividing $[M
_{\rm ab}\colon E]$. Therefore, if $M _{p} = E$, then this exponent is
not divisible by $p$. In view of the inclusion $N(M/E) \subseteq
N(R/E)$, $E ^{\ast }/N(R/E)$ is canonically isomorphic to a
homomorphic image of $E ^{\ast }/N(M/E)$, so the condition $M _{p} =
E$ ensures that the exponent $e(R/E)$ of $E ^{\ast }/N(R/E)$ is
also relatively prime to $p$. As $e(R/E)$ divides $[R\colon E]$, this
proves that $N(R/E) = E ^{\ast }$.
\par
Assume now that $R _{\rm ab} = E$ and $M _{p} \neq E$, denote by
$F _{1}$ the maximal abelian extension of $E$ in $M _{p}$, and by
$F _{2}$ the intermediate field of $M/E$ corresponding by Galois
theory to some Sylow $p$-subgroup of $G(M/E)$. Put $R _{1} = RF
_{1}$, $R _{2} = RF _{2}$ and $F _{3} = F _{1}F _{2}$. It follows
from Galois theory and the equality $R _{\rm ab} = E$ that the
compositum $RM _{p}$ is a Galois extension of $R$ with $G((RM
_{p})/R)$ canonically isomorphic to $G(M _{p}/E)$; in addition, it
becomes clear that $R _{1}$ is the maximal abelian extension of
$R$ in $RM _{p}$. Thus it turns out that $[R _{1}\colon R] = [F
_{1}\colon E]$, which means that $[R _{1}\colon E] = [R\colon 
E].[F _{1}\colon E]$. Observing that $[F _{2}\colon E]$ is not 
divisible by $p$, one also sees that $[R _{2}\colon F _{2}] = 
[R\colon E]$, $[(R _{1}F _{2})\colon F _{2}] = [R _{1}\colon E]$ 
and $[(RF _{3})\colon F _{2}] = [R _{2}\colon F _{2}].[F 
_{3}\colon F _{2}]$. The concluding equality and the normality of 
$F _{3}$ over $F _{2}$ imply that $R _{2} \cap F _{3} = F _{2}$. 
In view of Proposition 2.5 (iii) and Lemma 2.4, this leads to the 
conclusion that $N(R _{2}/E)N(F _{3}/E) = N(F _{2}/E)$. Note also 
that $N(F _{1}/E) = N(R _{1}/E)$. Indeed, it follows from Galois 
theory and the definition of $M _{p}$ that $M _{p}$ does not admit 
proper $p$-extensions in $M$, and by the inductive hypothesis, 
this yields $N((RM _{p})/M _{p}) = M _{p} ^{\ast }$. Hence, by the
former part of (1.1) (ii) and the transitivity of norm mappings,
we have $N((RM _{p})/E) = N(M _{p}/E) = N(F _{1}/E)$. At the same
time, since $R _{1}$ is the maximal abelian extension of $R$ in
$RM _{p}$, it turns out that $N((RM _{p})/R) = N(R _{1}/R)$,
which implies that $N((RM _{p})/E) = N(R _{1}/E) = N(F _{1}/E)$,
as claimed. The obtained results and the inclusions $N(R _{2}/E)
\subseteq N(R/E)$ and $N(F _{3}/E) \subseteq N(F _{1}/E)$,
indicate that $N(F _{2}/E)$ is a subgroup of $N(R/E)N(F _{1}/E) =$ 
\par \noindent
$N(R/E)N(R _{1}/E) = N(R/E)$. As $E ^{\ast }/N(R/E)$ and $E ^{\ast
}/N(F _{2}/E)$ are groups of finite relatively prime exponents,
this means that that $N(R/E) = E ^{\ast }$, so the proof of Lemma
4.1 is complete.
\par
\medskip
We are now in a position to prove Theorem 3.1 in the special case
where $R$ is an intermediate field of a finite Galois extension
$M/E$ with a solvable Galois group. It is clearly sufficient to
establish our assertion under the additional hypothesis that $N
_{p} (R _{1}/E _{1})$ and $N(R _{1}/E _{1})$ are related in
accordance with Theorem 3.1 whenever $E _{1}$ and $R _{1}$ are
extensions of $E$ in $R$ and $M$, respectively, such that $E _{1}
\neq E$ and $E _{1} \subseteq R _{1}$. Suppose that $R \neq E$,
put $\Phi = R _{{\rm ab},p}$, if $R _{\rm ab,p} \neq E$, and denote
by $\Phi $ some proper extension of $E$ in $R$ of primary degree,
otherwise (the existence of $\Phi $ in the latter case follows
from Galois theory and the well-known fact that maximal subgroups
of solvable finite groups are of primary indices). Also, let
$\alpha $ be an element of $N _{p} (R/E) \cap N(R _{{\rm
ab},p}/E)$, $\Phi ^{\prime }$ the maximal abelian $p$-extension of
$\Phi $ in $R$, $M ^{\prime }$ the compositum $\Phi M _{{\rm
ab},p}$, $k$ the maximal integer dividing $[M\colon E]$ and not
divisible by $p$, and $\Phi ^{\ast k} = \{z ^{k}\colon \ z \in \Phi
^{\ast }\}$. It is not difficult to see that $\Phi ^{\prime } \cap M
^{\prime } = \Phi $. Applying Proposition 2.5 (iii) or Lemma 2.4,
depending on whether or not Br$(\Phi ) _{p} \neq \{0\}$, one obtains
further that $\Phi ^{\ast } = N(\Phi ^{\prime }/\Phi )N(M ^{\prime
}/\Phi )$. Hence, by the inductive hypothesis and the inclusion
$N _{p} (M/\Phi ) \subseteq N _{p} (R/\Phi )$, $\Phi ^{\ast k}$ is
a subgroup of $N(R/\Phi ).(N(M ^{\prime }/\Phi ) \cap \Phi ^{\ast
k})$. Note also that Lemma 4.1 and the choice of $\Phi $ ensure
the existence of an element $\xi \in \Phi $ of norm $\alpha $ over
$E$. Taking now into account that $N(M ^{\prime }/E) \subseteq N(M
_{\rm ab,p}/E)$, one obtains that $\alpha ^{k} \in N(R/E)(N _{p}
(M/E) \cap N(M _{{\rm ab},p}/E))$, and then deduces from Lemma 3.3
that $\alpha ^{k} \in N(R/E)$. In view of the choice of $\alpha $
and $k$, this means that $\alpha \in N(R/E)$, which proves
Theorem 3.1 in the discussed special case. In order to do the
same in full generality, we need the following lemma.
\par
\medskip
{\bf Lemma 4.2.} {\it Let $E$ and $p$ satisfy the conditions of
Theorem 3.1, and let $R$ be an intermediate field of a finite
Galois extension $M/E$, such that $G(M/E) = [G(M/E),$ 
\par \noindent
$G(M/E)]$.
Then $N _{p} (R/E) \subseteq N(R/E)$.}
\par
\medskip
{\it Proof.} It is clearly sufficient to consider only the special
case of $R = M \neq E$ (and Br$(E) _{p} \neq \{0\}$). Denote by $E
_{p}$ be the intermediate field of $M/E$ corresponding by Galois
theory to some Sylow $p$-subgroup of $G(M/E)$. Then $p$ does not
divide the degree $[E _{p}\colon E] := m _{p}$, so the condition Br$(E)
_{p} \neq \{0\}$ guarantees that Br$(E _{p}) _{p} \neq \{0\}$. We
first show that $E ^{\ast } \subseteq N(M/E _{p})$, assuming
additionally that char$(E) = p$ or $E$ contains a primitive root
of unity of degree $[M\colon E _{p}]$. As $E$ is a quasilocal field, the
nontriviality of Br$(E _{p}) _{p}$ ensures that $E _{p}$ admits
local $p$-class field theory. Hence, by the former part of (1.1)
(ii), it is sufficient to prove the inclusion $E ^{\ast }
\subseteq N(L/E _{p})$, for an arbitrary cyclic extension $L$ of
$E _{p}$ in $M$. By [20, Sect. 15.1, Proposition b], this is
equivalent to the assertion that the cyclic $E _{p}$-algebra $(L/E
_{p}, \sigma , c)$ is isomorphic to the matrix $E _{p}$-algebra $M
_{n} (E _{p})$, where $c \in E ^{\ast }$, $n = [L\colon E _{p}]$ and
$\sigma $ is an $E _{p}$-automorphism of $L$ of order $n$. Since
g.c.d.$([E _{p}\colon E], p) = 1$, the surjectivity of the natural
homomorphism of Br$(E) _{p}$ into Br$(E _{p}) _{p}$ implies that
the corestriction homomorphism cor$_{E _{p}/E}\colon \ {\rm Br}(E 
_{p}) \to $ Br$(E)$ induces an isomorphism of Br$(E _{p}) _{p}$ on 
Br$(E) _{p}$ (cf. [27, Theorem 2.5]). Observe now that cor$_{E 
_{p}/E}$ maps the similarity class $[(L/E _{p}, \sigma , c)]$ into
$[(\widetilde L/E, \tilde \sigma , c)]$, for some cyclic
$p$-extension $\widetilde L$ of $E$ in $M$ (and a suitably chosen
generator $\tilde \sigma $ of $G(\widetilde L/E)$). Since $E$
contains a primitive root of unity of degree $[M\colon E _{p}]$ or
char$(E) = p$, this can be obtained by applying the projection
formula (cf. [16, Proposition 3 (i)] and [27, Theorem 3.2]), as
well as Kummer$^{,}$s theory and its analogue, due to Witt, for
finite abelian $p$-extensions over a field of characteristic $p$,
(see, for example, [13, Ch. 7, Sect. 3]). As $G(M/E) = [G(M/E),
G(M/E)]$, or equivalently, $M _{\rm ab} = E$, the obtained result
shows that $\widetilde L = E _{p}$ and $[(\widetilde L/E, \tilde
\sigma , c)] = 0$ in Br$(E)$. Furthermore, it becomes clear that
$[(L/E _{p}, \sigma , c)] = 0$ in Br$(E _{p})$, i.e. $c \in N(L/E
_{p})$, which proves the inclusion $E ^{\ast } \subseteq N(M/E
_{p})$. Since $N _{E} ^{E _{p}} (c) = c ^{m _{p}}$, one also sees
that $c ^{m _{p}} \in N(M/E)$, for each $c \in E ^{\ast }$.
\par
Suppose now that $p \neq $ char$(E)$, fix a primitive root of
unity $\varepsilon \in M _{\rm sep}$ of degree $[M\colon E _{p}]$, and
put $\Phi (\varepsilon ) = \Phi ^{\prime }$, for every
intermediate field $\Phi $ of $M/E$, and $H ^{m _{p}} = \{h ^{m
_{p}}\colon \ h \in H\}$, for each subgroup $H$ of $M ^{\prime \ast }$.
As $E ^{\prime }/E$ is an abelian extension, our assumption on
$G(M/E)$ ensures that $E ^{\prime } \cap M = E$, and by Galois
theory, this means that $M ^{\prime }/E ^{\prime }$ is a Galois
extension with $G(M ^{\prime }/E ^{\prime })$ canonically
isomorphic to $G(M/E)$. Thus it becomes clear from the previous
considerations that $E ^{\prime \ast m _{p}} \subseteq N(M
^{\prime }/E ^{\prime })$ and $N(E ^{\prime }/E) ^{m _{p}}
\subseteq N(M ^{\prime }/E) \subseteq N(M/E)$. Our argument also
shows that $M \cap E _{p} ^{\prime } = E _{p}$, and since $E _{p}$
is $p$-quasilocal, it enables one to deduce from Proposition 2.5
(iii), the former part of (1.1) (ii), and Lemma 2.2 that $N(M/E
_{p})N(E _{p} ^{\prime }/E _{p}) = E _{p} ^{\ast }$. Hence, by the
transitivity of norm mappings, $N(M/E)N(E _{p} ^{\prime }/E) = N(E
_{p}/E)$. These observations prove the inclusions $E ^{\ast m _{p}
^{2}} \subseteq N(E _{p}/E) ^{m _{p}} \subseteq N(M/E) ^{m
_{p}}.N(E ^{\prime }/E) ^{m _{p}} \subseteq N(M/E)$. This,
combined with the fact that $p$ does not divide $m _{p}$ and
the exponent of $E ^{\ast }/N(M/E)$ divides $[M\colon E]$, indicates
that $E ^{\ast m _{p}} \subseteq N(M/E)$ and so completes the
proof of Lemma 4.2.
\par
\medskip
It is now easy to accomplish the proof of Theorem 3.1. Assume that
$M _{0}$ is the maximal Galois extension of $E$ in $M$ with a
solvable Galois group, and also, that $\mu _{p}$, $m _{p}$ and
$\rho _{p}$ are the maximal integers not divisible by $p$
and dividing $[M _{0}\colon E]$, $[M\colon E]$ and $[R\colon E]$, 
respectively. Applying Lemma 3.3 to $M _{0}/E$ and Lemma 4.2 to 
$M/M _{0}$, one obtains that $E ^{\ast \mu _{p}} \subseteq N(M 
_{0}/E)$ and $M _{0} ^{\ast \bar m _{p}} \subseteq N(M/M _{0})$, 
where $\bar m _{p} = m _{p}/\mu _{p}$. Hence, by the norm identity 
$N _{E} ^{M} = N _{E} ^{M _{0}} \circ N _{M _{0}} ^{M}$, we have 
$E ^{\ast m _{p}} \subseteq N(M _{0}/E) ^{\bar m _{p}} \subseteq 
N(M/E)$. Since $E ^{\ast [R\colon E]} \subseteq N(R/E)$, $N(M/E) 
\subseteq N(R/E)$ and g.c.d.$(m _{p}, [R\colon E]) = \rho _{p}$, 
this means that $E ^{\ast \rho _{p}} \subseteq N(R/E)$, so Theorem 
3.1 is proved.
\par
\medskip
{\bf Remark 4.3.} Lemma 4.2 remains valid (with a slightly
modified proof), if the condition on $G(M/E)$ is replaced by the
one that $p$ does not divide the index $\vert G(M/E)\colon [G(M/E),
G(M/E)]\vert $. Note also that Lemma 3.3 can be deduced from (3.3)
and this generalization of Lemma 4.2, which allows us to skip
Lemma 3.2 and shorten the proof of Theorem 3.1. When $G _{E}$ is 
a prosolvable group, however, the inclusion of Lemma 3.2 enables 
us to deduce the theorem from Proposition 2.5 (iii), fundamentals 
of Galois theory, basic properties of cyclic algebras and 
well-known elementary facts concerning solvable finite groups. 
The prosolvability of $G _{E}$ is guaranteed, if $E$ possesses a 
Henselian discrete valuation (cf. [3, Corollary 2.5 and 
Proposition 3.1]).
\par
\medskip
{\it Proof of Theorem 1.1.} Since Theorem 1.1 (i) is a special case
of Theorem 3.1, it is sufficient to prove Theorem 1.1 (ii). Let $E$
be a quasilocal field and $R$, $\Phi (R)$ be finite extensions of
$E$ in $E _{\rm sep}$, such that $N(R/E) = N(\Phi (R)/E)$ and $\Phi
(R)/E$ is abelian. Applying Lemmas 2.2 and 2.3, one reduces the
proof of Theorem 1.1 (ii) to the special case in which Br$(E) _{p}
\neq \{0\}$, when $p$ ranges over the set $\Pi $ of prime numbers
dividing $[\Phi (R)\colon E]$. Let $\Lambda $ be the normal closure of $R$
in $E _{\rm sep}$ over $E$, and for each $p \in \Pi $, let $\Phi (R)
_{p}$ be the maximal $p$-extension of $E$ in $\Phi (R)$, $H _{p}$ be
a Sylow $p$-subgroup of $G(\Lambda /R)$, $G _{p}$ a Sylow
$p$-subgroup of $G(\Lambda /E)$ including $H _{p}$, $R _{1}$ and $E
_{1}$ the intermediate fields of $\Lambda /E$ corresponding by
Galois theory to $H _{p}$ and $G _{p}$, respectively. Note first
that $R _{\rm ab,p}$ is a subfield of $\Phi (R) _{p}$. Indeed, the
nontriviality of Br$(E) _{p}$ and the PQL-property of $E$ ensure the
availability of a local $p$-class field theory on $E$, so our
assertion follows from the fact that $N(\Phi (R)/E) = N(R/E)
\subseteq N(R _{\rm ab,p}/E)$ (whence, by Lemma 2.2, we have $N(\Phi
(R) _{p}/E) \subseteq N(R _{\rm ab,p}/E)$). It is easily verified
that $p$ does not divide $[R _{1}\colon R][E _{1}\colon E]$ and $R 
_{\rm ab,p}E _{1} = (\Phi (R) _{p}E _{1}) \cap R _{1}$. One also 
sees that Br$(E _{1}) _{p} \neq \{0\}$ (cf. [20, Sect. 13.4]). As 
$E$ is quasilocal, this indicates that $E _{1}$ admits local 
$p$-class field theory, so it follows from the former part of 
(1.1) (ii) that $N((R _{\rm ab,p}E _{1})/E _{1}) = N((\Phi (R) 
_{p}E _{1})/E _{1})N(R _{1}/E _{1})$. Our argument also proves 
that $N((R _{\rm ab,p}E _{1})/E) = N((\Phi (R) _{p}E _{1})/E)N(R 
_{1}/E) \subseteq (N(\Phi (R) _{p}/E)N(R/E) \cap N(E _{1}/E)) =$
\par \noindent
$N(\Phi (R) _{p}E _{1})/E)$. On the other hand, the inclusion $R 
_{\rm ab,p} \subseteq \Phi (R) _{p}$ implies that 
\par \noindent 
$N((\Phi (R) _{p}E _{1})/E) \subseteq N((R _{\rm ab,p}E _{1})/E)$, 
so it turns out that $N(\Phi (R) _{p}E _{1})/E) =$ \par \noindent 
$N((R _{\rm ab,p}E _{1})/E)$ and the quotient group $N(R _{\rm 
ab,p}/E)/N(\Phi (R) _{p}/E)$ is of exponent $e _{p}$ dividing $[E 
_{1}\colon E]$. Since $e _{p}$ divides $[\Phi (R)\colon E]$ and 
$p$ does not divide $[E _{1}\colon E]$, this means that $e _{p} = 
1$, i.e. $N(\Phi (R) _{p}/E) = N(R _{\rm ab,p}/E)$ (and $\Phi (R) 
_{p} = R _{\rm ab,p}$), for each $p \in \Pi $. Let now $p ^{\prime 
}$ be an arbitrary prime number. It is clear from the inclusion $R 
_{ab,p'} \subseteq R$ that $N(R/E) = N(\Phi (R)/E) \subseteq N(R 
_{ab,p'}/E)$ and $E ^{\ast }/N(R _{ab,p'}/E)$ is a homomorphic 
image of $E ^{\ast }/N(\Phi (R)/E)$, so it follows from Lemma 2.2 
that $E ^{\ast }/N(R _{ab,p'}/E)$ is a group of exponent dividing 
$[\Phi (R)\colon E]$. It is now easy to see that $N(R _{ab,p'}/E) 
= E ^{\ast }$ whenever $p ^{\prime } \not\in \Pi $, and to 
conclude that $N(R/E) = N(R _{ab}/E)$, as claimed by Theorem 1.1 
(ii).
\par
\medskip
{\bf Remark 4.4.} (i) The conditions of Theorem (i) are in force, 
if $E$ is a field with local class field theory in the sense of 
Neukirch-Perlis [19], i.e. if the triple $(G _{E}, \{G(E _{\rm 
sep}/F), F \in \Sigma \}, E _{\rm sep} ^{\ast })$ is an 
Artin-Tate class formation (cf. [2, Ch. XIV]), where $\Sigma $ is 
the set of finite extensions of $E$ in $E _{\rm sep}$. Then the 
assertion of Theorem 1.1 (i) is contained in [2, Ch. XIV, Theorem 
7]; in particular, it applies to any $p$-adically closed field 
and includes (1.1) (i) as a special case (see [21, Theorem 3.1 and 
Lemma 2.9] and [26, Ch. XIII, Proposition 6], respectively). 
\par
(ii) Let us note that the class of fields satisfying the 
conditions of Theorem 1.1 (i) is larger than the one studied in 
[19]. More precisely, for every divisible abelian torsion group 
$T$, there exists a quasilocal field $E(T)$ of this type, such 
that Br$(E(T))$ is isomorphic to $T$ and all finite groups are 
realizable as Galois groups over $E(T)$ (this will be proved 
elsewhere), whereas the Brauer groups of the fields considered in 
[19] embed in $\hbox{\Bbb Q}/{\Bbb Z}$. These properties of $E(T)$ 
indicate that it is strictly quasilocal if and only if the 
$p$-components of $T$ are nontrivial, for all prime numbers $p$.
\par
(iii) It follows at once from (1.1) (iii) and Theorem 1.1 (ii) that
$N(R/E)$ $= N(R _{\rm ab}/E)$ whenever $E$ is quasilocal and
algebraic over a global field $E _{0}$. In this case, $G _{E}$ is 
prosolvable and $E$ satisfies the conditions of Theorem 1.1 (i) 
as well (see [9, Proposition 2.7] and the references there).
\par
\vskip0.6truecm \centerline{\bf 5. Proof of Theorem 1.2}
\par
\medskip
In this Section we characterize (and prove the existence of)
Henselian discrete valued strictly quasilocal fields with the
properties required by Theorem 1.2. In what follows, $\overline 
P$ is the set of prime numbers, and for each field $E$, $P _{0} 
(E)$ is the subset of those $p \in \overline P$, for which $E$ 
contains a primitive $p$-th root of unity, or else, $p = $ 
char$(E)$. Also, we denote by $P _{1} (E)$ the subset of those $p
^{\prime } \in (\overline P \setminus P _{0} (E))$, for which $E
^{\ast } \neq E ^{\ast p'}$, and put $P _{2} (E) = \overline P
\setminus (P _{0} (E) \cup P _{1} (E))$. Every finite extension $L$
of a field $K$ with a Henselian valuation $v$ is considered with its
valuation extending $v$, this prolongation is also denoted by $v$
(unless stated otherwise), and $e(L/K)$ denotes the ramification
index of $L/K$. Our starting point is the following statement
(proved in [6]):
\par
\medskip
(5.1) With assumptions being as above, if $v$ is discrete, then
the following conditions are equivalent:
\par
(i) $K$ is strictly quasilocal;
\par
(ii) The residue field $\widehat K$ of $(K, v)$ is perfect, the
absolute Galois group $G _{\widehat K}$ is metabelian of
cohomological $p$-dimension cd$_{p}(G _{K}) = 1$, for each $p \in
\overline P$, and $P _{0} (\widetilde L) \subseteq P(\widetilde
L)$, for every finite extension $\widetilde L$ of $\widehat K$.
\par
\medskip
When these conditions are in force, $K$ is a nonreal field (cf.
[14, Theorem 3.16]), $P _{0} (K) \setminus \{{\rm char}(\widehat
K)\} = P _{0} (\widehat K) \setminus \{{\rm char}(\widehat K)\}$,
and the following is true:
\par
\medskip
(5.2) (i) Br$(\widetilde L) = \{0\}$ and Br$(L) _{p}$ is
isomorphic to the quasicyclic $p$-group $\hbox{\Bbb Z} (p
^{\infty })$, for every finite extension $L/K$ and each $p
\in P(\widehat K)$ (apply [25, Ch. II, Proposition 6 (b)] and
Scharlau's generalization of Witt's theorem [23]); in particular,
the natural homomorphism Br$(K) _{p} \to $ Br$(L) _{p}$ is
surjective;
\par
(ii) If $R$ is a finite extension of $K$ in $K _{\rm sep}$, such
that $[R\colon K]$ is not divisible by char$(\widehat K)$ or any $p
\in P _{2} (\widehat K)$, then $R$ is presentable as a compositum
of subextensions of $K$ of primary degrees; furthermore, if
$[R\colon K]$ is not divisible by char$(\widehat K)$ or any $p \in (P
_{1} (\widehat K) \cup P _{2} (\widehat K))$, then the normal
closure of $R$ in $K _{\rm sep}$ over $K$ has a nilpotent Galois
group (apply [10, (3.3)] and Galois theory);
\par
(iii) The group $G _{K}$ is pronilpotent in case char$(\widehat 
K) = 0$ and $P _{0} (\widehat K) = \overline P$.
\medskip
The following result (proved in [10]) sheds light on the norm
groups of finite separable extensions of a field $K$ subject to
the restrictions of (5.1). It shows that the conclusion of (1.1)
(i) is generally valid if and only if $P(\widehat K) = \overline
P$, i.e. $\widehat K$ is quasifinite.
\par
\medskip
{\bf Proposition 5.1.} {\it Assume that $(K, v)$ is a Henselian
discrete valued strictly quasilocal field, and $R$ is a finite
extension of $K$ in $K _{\rm sep}$. Then $R/K$ possesses an
intermediate field $R _{1}$ such that:}
\par
(i) {\it The sets of  prime divisors of $e(R _{1}/K)$, $[\widehat
R _{1}\colon \widehat K]$, $[\widehat R\colon \widehat R _{1}]$ 
and $[R\colon R _{1}]$ are included in $P _{1} (\widehat K)$, 
$\overline P \setminus P(\widehat K)$, $P(\widehat K)$ and $P _{0} 
(\widehat K) \cup P _{2} (\widehat K)$, respectively;}
\par
(ii) $N(R/K) = N((R _{\rm ab}R _{1})/K)$ {\it and $K ^{\ast
}/N(R/K)$ is isomorphic to the direct sum $G(R _{\rm ab}/K) \times
(K ^{\ast }/N(R _{1}/K))$;} \par (iii) $K ^{\ast }/N(R/K)$ {\it is
of order $[R _{\rm ab}\colon K][R _{1}\colon K] = [(R _{\rm ab}R 
_{1})\colon K]$.}
\par
\medskip
Our next result characterizes the fields singled out by Theorem 
1.2 (i)-(ii) in the class of strictly quasilocal fields with 
Henselian discrete valuations: 
\par
\medskip
{\bf Proposition 5.2.} {\it For a strictly quasilocal field $K$
with a Henselian discrete valuation $v$, the following conditions
are equivalent:}
\par
(i) $G _{K}$ {\it and the finite extensions of $K$ have the
properties required by Theorem 1.2;}
\par
(ii) char$(\widehat K) = 0${\it , $P _{0} (\widehat K) =
P(\widehat K) \neq \overline P$ and $P _{1} (\widehat K) =
\overline P \setminus P _{0} (\widehat K)$.
\par
When this occurs, every finite extension $R$ of $K$ in $K _{\rm
sep}$ is presentable as a compositum $R = R _{0}R _{1}$, where $R
_{1}$ is determined in accordance with Proposition 5.1 (i) and
(ii), and $R _{0}$ is an intermediate field of $R/K$ of degree $[R
_{0}\colon K] = [R\colon R _{1}]$. Moreover, the Galois group of 
the normal closure $\widetilde R$ of $R$ in $K _{\rm sep}$ over 
$K$ is nilpotent if and only if $R = R _{0}$.}
\par
\medskip
{\it Proof.} The implication (ii)$\to $(i) follows from
Proposition 5.1 and the fact that $R _{1}$ is defectless
over $K$ [28, Propositions 2.2 and 3.1]. The concluding
assertions of Proposition 5.2 are implied by (5.2) (ii), so we
assume further that condition (i) is in force. Let $\pi $ be a
generator of the maximal ideal of the valuation ring of $(K, v)$.
It is easily deduced from Proposition 5.1 that if $p \in P _{2}
(\widehat K)$ or $p = $ char$(\widehat K)$ and $p \not\in P _{0}
(K)$, then the root field, say $M _{\pi }$, of the binomial $X
^{p} - \pi $ satisfies the equality $N(M _{\pi }/K) = N(M _{\pi ,
{\rm ab}}/K)$. At the same time, it follows from (2.2) that $G(M
_{\pi }/K)$ is nonabelian and isomorphic to a semidirect product
of a group of order $p$ by a cyclic group of order dividing $p -
1$. This indicates that $G(M _{\pi }/K)$ is nonnilpotent. The
obtained results contradict condition (i), and thereby, prove
that $P _{0} (\widehat K) \cup P _{1} (\widehat K) = \overline
P$. Hence, by Galois theory, (1.1) (i) and condition (i),
$\widehat K$ is an infinite field. It remains to be seen that
char$(\widehat K) = 0$ and $P _{0} (\widehat K) \neq \overline
P$. Suppose that char$(\widehat K) = q > 0$ and $q \in P _{0}
(K)$. Then condition (i), statement (5.1) and the infinity of
$\widehat K$ imply the existence of a primitive $p$-th root of
unity in $\widehat K$, for at least one prime number $p \neq q$.
In addition, it becomes clear that there exists a cyclic
inertial extension $L _{p}$ of $K$ in $K _{\rm sep}$ of degree
$p$. Let $v _{p}$ be the valuation of $L _{p}$ extending $v$. It
is easily obtained from Galois theory (cf. [15, Ch. VIII, Theorem
20]) and the Henselian property of $v$ that $L _{p}$ has a normal
basis $B _{p}$ over $K$, such that $v _{p} (b) = 0$, for all $b
\in B _{p}$. Denote by $B _{p} ^{\prime }$ the polynomial set
$\{X ^{q} - X - b\pi ^{-1}\colon \ b \in B _{p}\}$, if char$(K) = q$,
and put $B _{p} ^{\prime } = \{X ^{q} - (1 + b\pi )\colon \ b \in B
_{p}\}$, in the mixed-characteristic case. It follows from the
Artin-Schreier theorem, Capelli's criterion (cf. [15, Ch. VIII,
Sect. 9]) and the Henselian property of $v _{p}$ that $B _{p}
^{\prime }$ consists of irreducible polynomials over $L _{p}$.
Furthermore, one obtains from Kummer's theory (and the assumption
that $q \in P _{0} (K)$) that the root field $L _{p} ^{\prime }$
of $B _{p} ^{\prime }$ over $L _{p}$ is a Galois extension of $K$
of degree $q ^{p}p$. It follows from the definition of $L _{p}
^{\prime }$ that the Sylow $q$-subgroup $G(L _{p} ^{\prime }/L
_{p})$ of $G(L _{p} ^{\prime }/K)$, is normal and elementary
abelian. At the same time, it is clear from the choice of $B
_{p}$ that $G(L _{p} ^{\prime }/L _{p})$ possesses maximal
subgroups that are not normal in $G(L _{p} ^{\prime }/E)$. These
properties of $G(L _{p} ^{\prime }/L _{p})$ indicate that $G(L
_{p} ^{\prime }/E)$ is nonnilpotent. On the other hand, since $q$
and $p$ lie in $P(\widehat K)$, Theorem 3.1 and the latter
assertion of (5.2) (i) show that $N(L _{p} ^{\prime }/K) = N(L
_{p,{\rm ab}} ^{\prime }/K)$. Thus the hypothesis that
char$(\widehat K) \neq 0$ leads to a contradiction with condition
(i), so the proof of Proposition 5.2 can be accomplished by
applying (5.2) (iii).
\par
\medskip
{\bf Corollary 5.3.} {\it Let $(K, v)$ be a Henselian discrete
valued field satisfying the conditions of Proposition 5.2,
$\varepsilon _{p}$ a primitive $p$-th root of unity in $K _{\rm
sep}$, for each $p \in \overline P$, and $[K(\varepsilon _{p})\colon K]
= \gamma _{p}$ in case $p \in (\overline P \setminus P(\widehat
K))$. Then each finite extension $L$ of $K$ in $K _{\rm sep}$ is
subject to the following alternative:}
\par
(i) $G _{L}$ {\it and finite extensions of $L$ have the properties
required by Theorem 1.2;} (ii) $G _{L}$ {\it is pronilpotent.
\par \noindent
The latter occurs if and only if the set $\Gamma (K) = \{\gamma
_{p}\colon \ p \in (\overline P \setminus P(\widehat K))\}$ is bounded
and $L$ contains as a subfield the inertial extension of $K$ in
$\overline K$ of degree equal to the least common multiple of the
elements of $\Gamma (K)$.}
\par
\medskip
{\it Proof.} Statement (5.1) and our assumptions guarantee that $P
_{0} (\widehat L) \cup P _{1} (\widehat L) = \overline P$, so the
stated alternative is contained in (5.2) (iii).
\par
\medskip
Our next result supplements Proposition 5.1 and combined with
Proposition 5.2, proves Theorem 1.2.
\par
\medskip
{\bf Proposition 5.4.} {\it Let $P _{0}$, $P _{1}$, $P _{2}$ and
$P$ be subsets of the set $\overline P$ of prime numbers, such
that $P _{0} \cup P _{1} \cup P _{2} = \overline P$, $2 \in P
_{0}$, $P _{i} \cap P _{j} = \phi \colon \ 0 \le i < j \le 2$, and $P
_{0} \subseteq P \subseteq (P _{0} \cup P _{2})$. For each $p \in
(P _{1} \cup P _{2})$, let $\gamma _{p}$ be an integer $\ge 2$
dividing $p - 1$ and not divisible by any element of $\overline P
\setminus P$. Assume also that $\gamma _{p} \ge 3$ in case $p \in
(P _{2} \setminus P)$. Then there exists a Henselian discrete
valued field $(K, v)$ satisfying the following conditions:}
\par
(i) $K$ {\it is strictly quasilocal with $P(\widehat K) = P$ and
$P _{j} (\widehat K) = P _{j}\colon \ j = 0, 1, 2$;}
\par
(ii) {\it For each $p \in (P _{1} \cup P _{2})$, $\gamma _{p}$
equals the degree $[K(\varepsilon _{p})\colon K]$, where $\varepsilon
_{p}$ is a primitive $p$-th root of unity in $K _{\rm sep}$.}
\par
\medskip
{\it Proof.} Denote by $G _{1}$ and $G _{0}$ the topological group
products $\prod _{p \in P} \hbox{\Bbb Z} _{p}$ and $\prod _{p \in
(\overline P \setminus P)} \hbox{\Bbb Z} _{p}$ (i.e. $G _{0} =
\{1\}$ in case $P = \overline P$), respectively, and fix an
algebraic closure $\overline {\hbox{\Bbb Q}}$ of the field of
rational numbers as well as a primitive $p$-th root of unity
$\varepsilon _{p} \in \overline {\hbox{\Bbb Q}}$, for each $p \in
\overline P$. Also, let $E _{0}$ be a subfield of $\overline
{\hbox{\Bbb Q}}$, such that $P _{0} (E _{0}) = P _{0}$, $P(E _{0})
= P$, $[E _{0} (\varepsilon _{p})\colon E _{0}] = \gamma 
_{p}\colon \ p \in (\overline P \setminus P _{0})$, and $G _{E 
_{0}} \cong G _{1}$ (the existence of $E _{0}$ is guaranteed by 
[6, Lemma 3.5]). Suppose further that $\varphi $ is a topological 
generator of $G _{E _{0}}$, and for each $p \in (\overline P 
\setminus P)$, $\delta _{p}$ is a primitive $\gamma _{p}$-th root 
of unity in $\hbox{\Bbb Z} _{p}$, $s _{p}$ and $t _{p}$ are integers, 
such that $\varphi (\varepsilon _{p}) = \varepsilon _{p} ^{s
_{p}}$, $t _{p} - \delta _{p} \in p\hbox{\Bbb Z} _{p}$, and $0 \le
s _{p}, t _{p} \le (p - 1)$. Assume also that the roots $\delta
_{p}$ are taken so that $t _{p} = s _{p}$ if and only if $p \in P
_{1}$. Regarding $\hbox{\Bbb Z} _{p}$ as as subgroup of $G _{0}$,
whenever $p \in (\overline P \setminus P)$, consider the
topological semidirect product $G = G _{0} \times G _{E _{0}}$,
defined by the rule $\varphi \lambda _{p}\varphi ^{-1} = \delta
_{p}\lambda _{p}\colon \ p \in (\overline P \setminus P)$, $\lambda
_{p} \in \hbox{\Bbb Z} _{p}$. It has been proved in [6, Sect. 3]
that there exists a Henselian discrete valued strictly quasilocal
field $(K, v)$, such that $G _{\widehat K}$ is continuously
isomorphic to $G$, $E _{0}$ is a subfield of $\widehat K$, and $E
_{0}$ is algebraically closed in $\widehat K$. In particular,
this implies that $P _{0} (\widehat K) = P _{0}$, $P(\widehat K) =
P$ and $[K(\varepsilon _{p})\colon K] = \gamma _{p}\colon \ p \in 
(\overline P \setminus P _{0})$. Applying finally (2.2) (ii), one 
concludes that $P _{1} (\widehat K) = P _{1}$ and so completes 
the proof of Proposition 5.4.
\par
\medskip
{\bf Corollary 5.5.} {\it There exists a set $\{(K _{n}, v 
_{n})\colon \ n \in \hbox{\Bbb N} \cup {\infty }\}$ of Henselian 
discrete valued strictly quasilocal fields satisfying the following 
conditions:}
\par
(i) {\it The absolute Galois group of a finite extension $R _{n}$
of $K _{n}$ is pronilpotent if and only if $n \in \hbox{\Bbb N}$
and $R _{n}$ contains as a subfield an inertial extension of $K
_{n}$ of degree $n$;}
\par
(ii) {\it Finite extensions of $K _{n}$ are subject to the
alternative described in Theorem 1.2, provided that $n \ge 2$.}
\par
\medskip
{\it Proof.} This follows at once from Corollary 5.3 and
Proposition 5.4.
\par
\medskip
{\bf Corollary 5.6.} {\it Let $P _{0}$ and $P$ be subsets of the
set $\overline P$ of prime numbers, such that $2 \in P _{0}$ and
$P _{0} \subseteq P$. Then there exists a strictly quasilocal
nonreal field $E$ such that:}
\par
(i) $P _{0} (E) = P _{0}$ {\it and $\{p \in \overline P\colon \ 
{\rm cd}_{p} (G _{E}) \neq 0\} = P$;}
\par
(ii) {\it If $P \neq P _{0}$, then $G _{E}$ is nonnilpotent and
finite extensions of $E$ are subject to the alternative described
by Theorem 1.2.}
\par
\medskip
{\it Proof.} Proposition 5.4 implies the existence of a Henselian
discrete valued strictly quasilocal field $(K, v)$, such that
char$(\widehat K) = 0$, $P _{0} (\widehat K) = P _{0}$, $P _{1}
(\widehat K) = \overline P \setminus P _{0}$, and for each $p \in
P _{1} (\widehat K)$, the extension of $K$ in $K _{\rm sep}$
obtained by adjoining a primitive $p$-th root of unity is of even
degree. By [3, Proposition 3.1], $G _{K}$ is a prosolvable group,
which means that it possesses a closed Hall pro-$P$-subgroup $H
_{P}$. Note finally that one can take as $E$ the intermediate 
field of $K _{\rm sep}/K$ corresponding by Galois theory to $H 
_{P}$.
\par
\vskip0.6truecm
\centerline{\bf Acknowledgements}
\par
\medskip
The proof of Theorem 1.1 was obtained during my visit to 
Tokai University, Hiratsuka, Japan (September 2002-March 2003). I 
gratefully acknowledge the stimulating atmosphere at the 
University as well as the ensured efficient fully supported 
medical care after my involvement into a travel accident. I would 
also like to thank my host professor M. Tanaka, Mrs M. Suzuki, Mrs 
A. Uchida and the colleagues from the Department of Mathematics 
for their kind hospitality. 
\vskip1cm \centerline{ REFERENCES} \vglue15pt\baselineskip12.8pt
\def\num#1{\smallskip\item{\hbox to\parindent{\enskip [#1]\hfill}}}
\parindent=1.38cm
\num{1} A. {\pc ALBERT}: {\sl Modern Higher Algebra.} Chicago Univ. Press, Chicago,
Ill., 1937.
\par
\num{2} E. {\pc ARTIN}; J. {\pc TATE}: {\sl Class Field Theory.} Benjamin, New
York-Amsterdam, 1968.
\par
\num{3} I.D. {\pc CHIPCHAKOV}: {\sl Henselian valued quasilocal fields with
totally indivisible value groups.} Comm. Algebra {\bf 27} (1999),
3093-3108.
\par
\num{4} I.D. {\pc CHIPCHAKOV}: {\sl Central division algebras of
$p$-primary dimensions and the $p$-component of the Brauer group
of a $p$-quasilocal field.} C.R. Acad. Sci. Bulg. {\bf 55} (2002),
55-60.
\par
\num{5} I.D. {\pc CHIPCHAKOV}: {\sl Norm groups and class fields of 
formally real quasilocal fields.} Preprint (available at  
www.arXiv.org/RA.math/0508019). 
\par
\num{6} I.D. {\pc CHIPCHAKOV}: {\sl Henselian discrete valued fields
admitting one-dimensional local class field theory.} In: 
Proceedings of AGAAP-Conference (V. Brinzanescu, V. Drensky and P.
Pragacz, Eds.), 23.9-02.10. 2003, Borovets, Bulgaria, Serdica
Math. J. {\bf 30} (2004), 363-394.
\par
\num{7} I.D. {\pc CHIPCHAKOV}: {\sl One-dimensional abstract local class
field theory.} Preprint (available at
www.math.arXiv.org/RA.math/0506515). 
\par
\num{8} I.D. {\pc CHIPCHAKOV}: {\sl On the residue fields of Henselian
valued stable fields.} Preprint (available at www.arXiv.org/RA.math/0412544). 
\par
\num{9} I.D. {\pc CHIPCHAKOV}: {\sl On nilpotent Galois groups and the
scope of the norm limitation theorem in one-dimensional abstract
local class field theory.} Preprint (to appear in the Proceedings 
of ICTAMI 05, held in Albac, Romania, 15.9-18.9, 2005; Acta 
Universitatis Apulensis {\bf 10}).
\par
\num{10} I.D. {\pc CHIPCHAKOV}: {\sl Class field theory for strictly
quasilocal fields with Henselian discrete valuations.} Preprint 
(available at www.arXiv.org/ \par RA.math/0506069). 
\par
\num{11} K. {\pc IWASAWA}: {\sl Local Class Field Theory.} Iwanami Shoten,
Japan, 1980 (Japanese: Russian transl. in Mir, Moscow, 1983). 
\par
\num{12} B. {\pc JACOB}; A. {\pc WADSWORTH}: {\sl  Division  algebras  over
Henselian fields.} J. Algebra {\bf 128} (1990), 126-179.
\par
\num{13} G. {\pc KARPILOVSKY}: {\sl Topics in Field Theory.} North-Holland
Math. Studies, 155, North Holland, Amsterdam etc., 1989.
\par
\num{14} T.Y. {\pc LAM}: {\sl Orderings, valuations and quadratic forms.}
Conf. Board Math. Sci. Regional Conf. Ser. Math. No 52, Am. Math.
Soc., Providence, RI, 1983.
\par
\num{15} S. {\pc LANG}: {\sl Algebra.} Addison-Wesley Publ. Comp., Mass.,
1965.
\par
\num{16} P. {\pc MAMMONE}; A. {\pc MERKURJEV}: {\sl On the corestriction of the
$p ^{n}$-symbol.} Comm. Algebra {\bf 76} (1991), 73-80.
\par
\num{17} G.A. {\pc MILLER}; H.C. {\pc MORENO}: {\sl Nonabelian groups in which
every subgroup is abelian.} Trans. Amer. Math. Soc. {\bf 4} (1903),
398-404.
\par
\num{18} M. {\pc MORIYA}: {\sl Eine notwendige Bedingung f$\ddot 
u$r die G$\ddot u$ltigkeit der \par Klassenk$\ddot o$rpertheorie 
im Kleinen.} Math. J. Okayama Univ. {\bf 2} (1952), 13-20.
\par
\num{19} J. {\pc NEUKIRCH}; R. {\pc PERLIS}: {\sl Fields with local class field
theory.} J. Algebra {\bf 42} (1976), 531-536.
\par
\num{20} R. {\pc PIERCE}: {\sl Associative Algebras.} Graduate Texts in 
Math. {\bf 88}, Springer-Verlag, New York-Heidelberg-Berlin, 1982.
\par
\num{21} A. {\pc PRESTEL}; P. {\pc ROQUETTE}: {\sl Formally $p$-adic Fields.}
Lecture Notes in Math. {\bf 1050}, Springer-Verlag, Berlin etc.,
1984.
\par
\num{22} L. {\pc REDEI}: {\sl Algebra, v. 1.} Akademiai Kiado, Budapest,
1967.
\par
\num{23} W. {\pc SCHARLAU}: {\sl $\ddot U$ber die Brauer-Gruppe eines
Hensel-K$\ddot o$rpers.} Abh. Math. Semin. Univ. Hamb. {\bf 33}
(1969), 243-249.
\par
\num{24} O.F.G. {\pc SCHILLING}: {\sl Necessary conditions for local class
field theory.} Math. J. Okayama Univ. {\bf 3} (1953), 5-10.
\par
\num{25} J.-P. {\pc SERRE}: {\sl Cohomologie Galoisienne.} Lecture Notes in
Math. {\bf 5}, Springer-Verlag, Berlin-Heidelberg-New York, 1965.
\par
\num{26} J.-P. {\pc SERRE}: {\sl Local Fields.} Graduate Texts in
Mathematics {\bf 67}, Springer-Verlag, New York-Heidelberg-Berlin,
1979.
\par
\num{27} J.-P. {\pc TIGNOL}: {\sl On the corestriction of central simple
algebras.} Math. Z. {\bf 194} (1987), 267-274.
\par
\num{28} I.L. {\pc TOMCHIN}; V.I. {\pc YANCHEVSKIJ}: {\sl On defects of valued
division algebras.} Algebra i Analiz {\bf 3} (1991), No 3, 147-164
(in Russian: Engl. transl. in St. Petersburg Math. J. {\bf 3}
(1992), No 3, 631-646).
\par
\num{29} G. {\pc WHAPLES}: {\sl Generalized local class field theory. II.
Existence Theorem. III. Second form of the existence theorem.
Structure of analytic groups. IV. Cardinalities.} Duke Math. J.
{\bf 21} (1954), 247-255, 575-581, 583-586.
\vskip1cm
\def\pc#1{\eightrm#1\sixrm}
\hfill\vtop{\eightrm\hbox to 5cm{\hfill Ivan {\pc
CHIPCHAKOV}\hfill}
 \hbox to 5cm{\hfill Institute of Mathematics and Informatics\hfill}\vskip-2pt
 \hbox to 5cm{\hfill Bulgarian Academy of Sciences\hfill}
\hbox to 5cm{\hfill Acad. G. Bonchev Str., bl. 8\hfill} \hbox to
5cm{\hfill 1113 {\pc SOFIA,} Bulgaria\hfill}}
\end